\documentclass[12pt]{amsart}

\usepackage{amsfonts}
\usepackage{amsmath}
\usepackage{amssymb}

\usepackage{mathrsfs}
\usepackage{amscd}
\usepackage[all]{xy}
\usepackage[pdftex,final]{graphicx}

\usepackage{hyperref}

\usepackage{exscale,txfonts}

\newtheorem{lem}{Lemma}[section]
\newtheorem{thm}[lem]{Theorem}
\newtheorem{prop}[lem]{Proposition}
\newtheorem{cor}[lem]{Corollary}
\theoremstyle{definition}

\newtheorem{rem}[lem]{Remark}

\usepackage{mathtools}

\makeatletter
\DeclareRobustCommand\widecheck[1]{{\mathpalette\@widecheck{#1}}}
\def\@widecheck#1#2{%
    \setbox\z@\hbox{\m@th$#1#2$}%
    \setbox\tw@\hbox{\m@th$#1%
       \widehat{%
          \vrule\@width\z@\@height\ht\z@
          \vrule\@height\z@\@width\wd\z@}$}%
    \dp\tw@-\ht\z@
    \@tempdima\ht\z@ \advance\@tempdima2\ht\tw@ \divide\@tempdima\thr@@
    \setbox\tw@\hbox{%
       \raise\@tempdima\hbox{\scalebox{1}[-1]{\lower\@tempdima\box
\tw@}}}%
    {\ooalign{\box\tw@ \cr \box\z@}}}
\makeatother
\newcommand{\Q}{\mathbb{Q}}

\newcommand{\Z}{\mathbb{Z}}

\newcommand{\C}{\mathbb{C}}

\newcommand{\spl}[2]{\mathrm{SL}_{#1}(#2)}

\newcommand{\gl}[2]{\mathrm{GL}_{#1}(#2)}

\newcommand{\gpb}[1]{\left[ #1\right]} 

\newcommand{\spec}[1]{S_{#1}}
\newcommand{\pn}[2]{\mathbb{P}^{#1}(#2)}
\newcommand{\projl}[1]{\pn{1}{#1}}

\newcommand{\<}{\langle}
\renewcommand{\>}{\rangle}

\renewcommand{\ker}[1]{\mathrm{Ker}(#1)}

\newcommand{\ntors}[2]{{#1}\left[ #2 \right]}

\renewcommand{\spec}[1]{\mathrm{Spec}(#1)}

\newcommand{\kind}[1]{K^{\mathrm{\small ind}}_3(#1)}

\newcommand{\ch}[2]{\bar{c}_{#1,#2}}

\newcommand{\ho}[3]{\mathrm{H}_{#1}(#2,#3 )}

\newcommand{\pont}{\ast} 

\title{The Chern class for $K_3$ and the cyclic quantum dilogarithm}
\author{Kevin Hutchinson}
\address{School of Mathematics and Statistics,
 University College Dublin}
\email{kevin.hutchinson@ucd.ie}
\date{\today}

\keywords{
$K$-theory, Bloch group, Chern class}
\subjclass{19F99, 20G10}

\begin{document}
\bibliographystyle{plain}
\maketitle
{

\begin{abstract}
In this note we confirm the conjecture of Calegari, Garoufalidis and Zagier in \cite{cgz:bloch} that $R_\zeta=c_\zeta^2$ where $R_\zeta$ is their map on $K_3$ defined using the cyclic quantum dilogarithm and $c_\zeta$ is  the Chern class map on $K_3$. 
\end{abstract}
\section{Introduction}
Let $N\geq 1$ be odd and let $\zeta$ be a primitive $N$th root of $1$. In a recent article  (\cite{cgz:bloch}), Calegari, Garoufalidfis and Zagier use the cyclic quantum dilogarithm and the relationship between the $K_3$ of a field and the Bloch group of the field to define a natural homomorphism $R_\zeta:K_3(F)\to F_N^\times/(F^\times_N)^N$ where  $F$ is a field 
such that $\mu_N(F)=\{ 1\}$ and $F_N:=F(\mu_N)$. The $K$-theory Chern class also provides such a homomorphism (even without the restriction on $\mu_F$),  which we will denote $c_\zeta$, following \cite{cgz:bloch}.  Calegari, Garoufalidis and Zagier show that these two  maps are closely related. To be precise, they show  that there exists $\gamma\in (\Z/N)^\times$, independent of the field $F$, such that $R_\zeta=c_\zeta^{\gamma}$  (\cite[Theorem 1.6]{cgz:bloch}). Furthermore,  in the introduction to the paper they conjecture that $\gamma =2$.  In this  note, we prove this conjecture. 

To do this, we will broadly follow the plan  sketched out by Calegari, Garoufalidis and Zagier in section 5.3 of their article. Using the good functorial properties of  the maps $R_\zeta$ and $c_\zeta$
they reduce  the calculation of $\gamma$ to the calculation of $c_\zeta(\eta_\zeta)$ where $\eta_\zeta$ is an explicitly given element of\\
 $K_3(\Q(\zeta+\zeta^{-1})/N$  (see section 3 below) when $N$ is a power of an odd prime $\ell$. One of the main results of their paper is that $R_\zeta(\eta_\zeta)=\zeta^2\in  F_N^\times/(F_N^\times)^N$ (\cite[Theorem 7.4]{cgz:bloch}).  Our Theorem \ref{thm:main} below shows that $c_\zeta(\eta_\zeta)=\zeta$, from which it follows 
(\cite[section 1.2 and  section 5.3]{cgz:bloch}) that $R_\zeta=c_\zeta^2$.

In section \ref{sec:review} we review some of the relevant facts about $K$-theory, homology of linear groups and Chern classes that we will require. In section \ref{sec:eta}, we review the definition and properties of the element $\eta_\zeta$. In section \ref{sec:bott} we show that (as guessed by Calegari, Garoufalidis and Zagier) the image of $\eta_\zeta$ in $K_3(\Q(\zeta);\Z/N)$ is $\zeta*\beta$ where $\zeta\in K_1(F)=F^\times$ and $\beta\in K_2(F;\Z/N)$ is the Bott element in $K$-theory. The result then follows from an application of a multiplication formula for Chern classes due to Soul\'e (see \ref{lem:sou}). 

We do not define the map $R_\zeta$ below, as it is not needed for the purposes of this article. For its definition and properties, see the article of Calegari, Garoufalidis and Zagier (\cite{cgz:bloch}). Their motivation for defining this map is because its values coincide with  certain algebraic numbers arising from the quantum modularity conjecture for the Kashaev 
invariant of knots. Furthermore, they use its properties to prove a conjecture of W. Nahm that the modularity of certain hypergeometric $q$-series implies that an associated element of the Bloch group of $\bar{\Q}$ vanishes. 


\section{Review of some background}\label{sec:review}

\subsection{The homology of a cyclic group}
Let $C$ be a cyclic group of order $N$ with generator $t$.

Let $N>1$ be an odd integer.  Let $B_\bullet(C)$ be the (right) bar resolution of $\Z$ over $\Z C$. 
 So $B_r(C)$ is a free right $\Z C$-module with basis $\{ [c_1|\cdots |c_r]\ |\ c_i\in C\}$ with $d_r:B_r(C)\to B_{r-1}(C)$ given by 
\[
d_r([c_1|\cdots |c_r]):= [c_2|\cdots |c_r]-[c_1c_2|\cdots |c_r]+\cdots +(-1)^{r-1}[c_1|\cdots |c_{r-1}c_r]+(-1)^r[c_1|\cdots c_{r-1}]c_r.
\]
For $r\geq 0$, let $s=\lfloor r/2\rfloor$ and let 
\[
\alpha_r(t):=
\left\{
\begin{array}{ll}
\sum_{j_1,\ldots,j_s=0}^{N-1}[t|t^{j_1}|\cdots |t|t^{j_s}],& r\mbox{ even}\\
\sum_{j_1,\ldots,j_s=0}^{N-1}[t|t^{j_1}|\cdots |t|t^{j_s}|t], & r\mbox{ odd}\\
\end{array}
\right.
\]
in $B_r(C)$. 

\begin{lem} \label{lem:cyc}  Let $C$ be cyclic of order $N$ with generator $t$.
\begin{enumerate}
\item Let $(S_r(C),d_r^{(t)})$ be the resolution where $S_r(C)=\Z C$ for all $r$ and  $d_r=d_r^{(t)}:S_r(C)\to S_{r-1}(C)$ is the 
map 
\[
d_r(1):=
\left\{
\begin{array}{ll}
1+t+\cdots +t^{N-1},& r \mbox{ even}\\
1-t,& r \mbox{ odd}.
\end{array}
\right.
\]
Then the maps 
\[
A_r:S_r(C)\to B_r(C), \quad 1\mapsto \alpha_r(t)
\]
define an augmentation preserving map of $\Z C$-complexes.
\item For all odd $r$ the class $[\alpha_r(t)]$ of $\alpha_r(t)\otimes 1 \in B_r(C)\otimes_{\Z C}\Z$ represents a generator of the 
cyclic group $H_r(C,\Z)$ of order $N$.
\item For all $r\geq 1$, the class of $\alpha_r(t)\otimes 1 \in B_r(C)\otimes_{\Z C}\Z/N$ represents a generator of the 
cyclic group $H_r(C,\Z/N)$ of order $N$.
\end{enumerate}
\end{lem}
\begin{proof} 
\begin{enumerate}
\item One easily verifies that 
\[
d_r(\alpha_r(t))=
\left\{
\begin{array}{ll}
\alpha_{r-1}(t)-\alpha_{r-1}(t)t,& r \mbox{  odd}\\
\sum_{j=0}^{N-1}\alpha_{r-1}(t)t^j,& r \mbox{ even}.\\
\end{array}
\right.
\]
\item When $r$ is odd, $1\in S_r(C)=\Z C$ represents a generator of\\
 $\ho{r}{C}{\Z}=H_r(S_\bullet(C)\otimes_{\Z C}\Z)$.
\item Likewise, when $r$ is even, $1\in S_r(C)\otimes_{\Z C} \Z/N$ represents a generator of\\
 $\ho{r}{C}{\Z/N}=H_r(S_\bullet(C)\otimes_{\Z C}\Z/N)$.
\end{enumerate}
\end{proof}

\begin{cor} \label{cor:beta} Let  $\tilde{\beta}=\tilde{\beta}(t)$ denote the generator $[\alpha_2(t)]\in H_2(C,\Z/N\Z)$. 
 Let $\pont$ denote the Pontryagin product on $H_\bullet(C,\Z/N\Z)$. Then $[t]\pont \tilde{\beta}=[\alpha_3(t)]$. 
\end{cor}

\begin{proof}
The Pontryagin product is described on the level of chains in the bar resolution by the shuffle product:
\[
[g_1|\cdots|g_p]\pont [g_{p+1}|\cdots|g_{p+q}]=\sum_{\sigma} (-1)^{\mbox{\small sgn($\sigma$)}}[g_{\sigma^{-1}(1)}|\cdots |g_{\sigma^{-1}(p+q)}]
\]
where $\sigma$ runs over all $(p,q)$-shuffles (see \cite[V.5]{brown:coh}). Thus 
\[
[\alpha_1(t)]\pont [\alpha_2(t)]=\sum_{j=0}^{N-1}\left( [t|t|t^j]-[t|t|t^j]+[t|t^j|t]\right)=\sum_{j=0}^{N-1}[t|t^j|t]=[\alpha_3(t)]
\]
as required. 
\end{proof}


\subsection{$K$-theory, the Hurewicz homomorphism and the Bloch group.}
We begin by reviewing some results about  the structure of $K_3(F)$ for number fields $F$.
Recall that for a field $E$, $w_2(E)$ denotes $\mathrm{sup}\{ m\geq 1 \ |\ E(\mu_m)/E \mbox{ is a multiquadratic extension}\}$. 
\begin{lem}\label{lem:k3} Let $F$ be a number field.
\begin{enumerate}
\item If $F$ is totally imaginary, then  ${K_3(F)}=\kind{F}\cong \Z/w_2(F)\oplus \Z^{s}$ where $s=[F:\Q]/2$. 
\item If $F$ has $r>1$ real embeddings and $s$ complex embeddings then 
\[
K_3(F)\cong \Z^s\oplus \Z/2w_2(F)\oplus (\Z/2)^{r-1}.
\]
\end{enumerate}
\end{lem}
\begin{proof}\ 
\begin{enumerate}
\item See \cite[VI.5.3]{weibel:kbook}.
\item See \cite[VI.5.3]{weibel:kbook}.
\end{enumerate}
\end{proof}
For a ring $R$ there are natural Hurewicz homomorphisms $h_n:K_n(R)\to \ho{n}{\gl{}{R}}{\Z}$.  Here $\gl{}{R}:=\lim_{n\to \infty}\gl{n}{R}$, the limit being taken with respect to the 
injective maps\\
 \[
j_n:\gl{n}{R}\to \gl{n+1}{R}, A\mapsto 
\left[
\begin{array}{cc}
A&0\\
0&1\\
\end{array}
\right].
\]

We summarize some results of Suslin relevant to our calculation below. Recall that the $K$-theory $K_\bullet(R)$ of a commutative ring $R$ is a graded commutative ring and that for a field $F$, there is an isomorphism $F^\times\cong K_1(F), a\mapsto \{ a\}$ (or $\ell(a)$). 

\begin{prop}\label{prop:sus} Let $F$ be an infinite field. Let $\bar{K}_3(F)$ denote $K_3(F)/(\{ -1\}\cdot K_2(F))$.
\begin{enumerate}
\item  The Hurewicz homomorphism $h_3: K_3(F)\to H_3(\gl{}{F},\Z)$ induces an isomorphism $\bar{h}_3:\bar{K}_3(F)\cong H_3(\spl{}{F},\Z)$. 
\item The inclusion $j:\gl{3}{F}\to \gl{}{F}$ induces an isomorphism $\ho{3}{\gl{3}{F}}{\Z}\cong \ho{3}{\gl{}{F}}{\Z}$.
\end{enumerate}
\end{prop}
\begin{proof}
\begin{enumerate}
\item This is a special case of  \cite[Corollary 5.2]{sus:bloch}.
\item This is a special case of the main  homology stability theorem of \cite{sus:homgln}.
\end{enumerate}
\end{proof}
\begin{rem} 
Note that $K_3(F)$ maps onto $\bar{K}_3(F)$ with kernel annihilated by $2$. In particular, for odd $n$ we have $K_3(F)/n=\bar{K}_3(F)/n$.
\end{rem}

\begin{rem}\label{rem:kind}  By Matsumoto's theorem, for any field $F$, the group $K_2(F)$ is generated by all products $\{ a \}\cdot \{ b\}:= \{ a,b\}$, with $a,b\in F^\times$. Furthermore, for any field $F$ there is short exact sequence $0\to K_3^M(F)\to K_3(F)\to \kind{F}\to 0$ where $K_3^M(F)$ is the subgroup generated by products 
$\{ a,b,c\}=\{ a\}\cdot \{ b\}\cdot \{ c\}$ for $a,b,c\in F^\times$. 

For a number field $F$, it is known (see \cite[p 400]{bass:tate}) that $K_3^M(F)$ is $2$-torsion generated by products $\{ -1,-1,c\}$. Hence $K_3^M(F)=
\{ -1\}\cdot K_2(F)$ and $\bar{K}_3(F)=\kind{F}$ when $F$ is a number field.
\end{rem}

The \emph{ pre-Bloch group} or \emph{scissors congruence group}  of the field $F$ is the additive abelian group $P(F)$ with generators $\gpb{x}$, $a\in F\setminus \{ 0,1\}$, subject to the relations 
\[
\gpb{x}-\gpb{y}+\gpb{\frac{y}{x}}-\gpb{\frac{1-x^{-1}}{1-y^{-1}}}+\gpb{\frac{1-x}{1-y}},\quad x\not=y \in F\setminus\{ 0,1\}.
\]
The \emph{Bloch group} of $F$ is the subgroup $B(F):=\ker{\lambda}$ where 
\[
\lambda:P(F)\to \frac{F^\times \otimes F^\times}{\langle x\otimes y+ y\otimes x | x,y \in F^\times \rangle}, \gpb{x}\mapsto x \otimes (1-x).
\]

\begin{rem}
 This is Suslin's definition of the Bloch group. As Calegari, Garoufalidis and Zagier point out, it differs by some $2$-torsion with the version used in their article. Since our main result concerns $B(F)/n$ where $n$ is odd, this difference does not concern us here. 
\end{rem}

For any infinite field $F$, there is a natural homorphism $\ho{3}{\gl{2}{F}}{\Z}\to B(F)$ first described by Dupont and Sah (\cite{sah:dupont}). In \cite{sus:bloch} (see section 3), Suslin 
constructs a homomorphism $\partial: \ho{3}{\gl{3}{F}}{\Z}\to B(F)$  which extends the homomorphism of Dupont and Sah. In view of  Proposition \ref{prop:sus} (2) above, this gives rise to composite homomorphism $\partial\circ h_3:K_3(F)\to B(F)$. 

For a finite cyclic group $C$ of even order, let $\tilde{C}$ denote the non-trivial extension by $\Z/2$. The following is \cite[Theorem 5.2]{sus:bloch}:
\begin{thm}\label{thm:suslin}  For any infinite field $F$, there is a natural short exact sequence 
\[
0\to \widetilde{\mu_F}\to \kind{F}\to B(F)\to 0.
\]
\end{thm}

\begin{cor}\label{cor:suslin} Let $F$ be an infinite field and let $N\geq 1$ be odd with $\mu_N(F):=\mu_N\cap \mu_F=\{ 1\}$. Then the map $K_3(F)\to B(F)$ induces an isomorphism
$\kind{F}/N\cong B(F)/N$. 
\end{cor}
\subsection{$K$-theory with finite coefficients and Chern classes} 
For $\ell\geq 1$, the $K$-theory groups  with finite coefficients $K_m(R;\Z/ \ell)$ of a ring $R$  fit functorially into a universal coefficient sequence
\[
0\to K_m(R)/\ell\to K_m(R;\Z/\ell)\to \ntors{K_{m-1}(R)}{\ell}\to 0
\]
where $\ntors{A}{\ell}:= \{ a\in A\ |\ \ell a=0\}$. (See for example \cite[IV.2]{weibel:kbook}.)  Furthermore,  for a commutative ring $R$, there is a natural product  $K_n(R;Z/N)\times K_m(R;\Z/N)\to K_{n+m}(R;\Z/N)$ making $K_\bullet(R;\Z/N)$ into a graded commutative ring (see \cite[II.2]{soule:etale}). 

Let $R$ be a commutative ring in which $N$ is invertible. For $i,k\geq 0$ let
\[
\ch{i}{k}:K_{2i-k}(R,\Z/N)\to H^k(\spec{R},\mu_N^{\otimes i}):=H^k(R,\mu_N^{\otimes i}).
\]
denote the $K$-theory \'etale Chern classes  (\cite[II.2.3]{soule:etale}).

For a field $F$, the Hurewicz homomorphism induces an isomorphism $K_1(F)\cong H_1(GL(F),\Z)\cong F^\times$ (the second isomorphism is induced by the determinant). Similarly, we obtain an isomorphism
$K_1(F,\Z/N)\cong H_1(GL(F),\Z/N)\cong F^\times/(F^\times)^N$. With this identification, the map $\ch{1}{1}:F^\times/(F^\times)^N=K_1(F,\Z/N)\to H^1(F,\mu_N)$ is just the Kummer isomorphism.

The following is \cite[Theorem 4:12]{levine:k3}
\begin{thm}\label{lem:levine}  For any field $E$, the Chern class $\bar{c}_{2,1}$  induces an isomorphism\\
 $K_3(E;\Z/n)^{\mbox{\tiny ind}}\cong H^1(E,\mu_n^{\otimes 2})$, provided $(n,\mathrm{char}(E))=1$.
\end{thm}

\begin{cor}\label{cor:levine}
Let $E$ be a  number field containing $\mu_N$. In this case, the Chern class 
$\ch{2}{1}:K_3(E,\Z/N)\to E^\times/(E^\times)^N\otimes \mu_N$ is an isomorphism. 
\end{cor}
\begin{proof}  Since $E$ is totally imaginary, $K_3^M(E)=0$ and hence $K_3(E;\Z/N)=K_3(E;\Z/N)^{\mbox{\tiny ind}}$.
\end{proof}

As in \cite{cgz:bloch}, let $F$ be a number field not containing a primitive $N$-th root of unity and let $F_N:=F(\mu_N)$.  Fix a primitive $N$th root of unity $\zeta=\zeta_N\in \mu_N$. 
Then the Chern class $\ch{2}{1}$ gives rise to a homomorphism $c_\zeta$
\[
\xymatrix{
K_3(F;\Z/N)\ar^-{\ch{2}{1}}[r]&H^1(F,\mu_N^{\otimes 2})\ar[r]& H^1(F_N,\mu_N^{\otimes 2})\cong (F_N^\times)/(F_N^\times)^N\otimes \mu_N\ar^-{\cong}[r]& (F_N^\times)/(F_N^\times)^N,
}
\]
where the last arrow sends $a\otimes \zeta$ to $a$.  Combining with the inclusion $K_3(F)/N\to K_3(F;\Z/N)$ gives the map 
\[
c_\zeta:K_3(F)/N\to (F_N^\times)/(F_N^\times)^N.
\]

We will require below the following multiplication formula for $K$-theory Chern classes due to Soul\'e:

\begin{thm}\label{thm:sou}(\cite[II.3, Th\'eor\`eme 1(i)]{soule:etale})
Let $q=\ell^m$ with $\ell$ prime. Let $A$ be a commutative ring in which $\ell$ is invertible. Let $a\in K_n(A;\Z/q)$, $b\in K_r(A;\Z/q)$, $n,r\geq 1$. Then
\[
\ch{i}{k}(a*b)=\sum -\frac{(i-1)!}{(i'-1)!(i''-1)!}\ch{i'}{k'}(a)\cup\ch{i''}{k''}(b)
\]
where the sum ranges over all $i',i'',k',k''$ satisfying $i=i'+i''$, $k=k'+k''$, $2i'-k'=n$, $2i''-k''=r$.
\end{thm}

\begin{cor}\label{cor:sou}
If $a\in K_1(A;\Z/q)$ and $b\in K_2(A;\Z/q)$ then $\ch{2}{1}(a*b)=-\ch{1}{1}(a)\cup\ch{1}{0}(b)$.
\end{cor}
\begin{proof} The conditions on $i',i'',k',k''$ imply that the only possible  values are $i'=1=i''$.
\end{proof}
\section{The element $\eta_\zeta$}\label{sec:eta}

Now let $N=\ell^m$, where $\ell$ is an odd prime number.   Fix a generator $\zeta$ of the group $\mu_N$ of $N$th roots of unity in $\C$.  Let $\Q(\zeta+\zeta^{-1}):=\Q(\zeta)^+$.  Observe that 
\[
w_2(\Q(\zeta)^+)=
\left\{
\begin{array}{ll}
24N,& \ell\not= 3,\\
8N,& \ell = 3.\\
\end{array}
\right.
\]
Thus
\[
K_3(\Q(\zeta)^+)/N\cong\bar{K}_3(\Q(\zeta)^+)/N=\kind{\Q(\zeta)^+}/N\cong B(\Q(\zeta)^+)/N\cong \Z/N
\]
by Lemma \ref{lem:k3} and Corollary \ref{cor:suslin}.

Let $t=t_\zeta:= \left[
\begin{array}{cc}
0&1\\
-1&\zeta+\zeta^{-1}
\end{array}
\right]\in \spl{2}{\Q(\zeta)^+}$.  In $\spl{2}{\Q(\zeta)}$, we have 
\[
t_\zeta= \left[
\begin{array}{cc}
\zeta&\zeta ^{-1}\\
1&1\\
\end{array}
\right]
 \left[
\begin{array}{cc}
\zeta ^{-1}&0\\
0&\zeta\\
\end{array}
\right]
 \left[
\begin{array}{cc}
\zeta&\zeta ^{-1}\\
1&1\\
\end{array}
\right]^{-1}.
\]
In  particular, $t\in \spl{2}{\Q(\zeta)^+}$ is an element of order $N$.

Now  we define  $\eta_\zeta$    the image of the generator $[\alpha_3(t)]\in \ho{3}{\langle t\rangle}{\Z}$ under the composite map 
\[
\xymatrix{
\ho{3}{\langle t\rangle}{\Z}\ar[r]& \ho{3}{\spl{}{\Q(\zeta)^+}}{\Z}\ar^-{\bar{h}_3^{-1}}[r]&\bar{K}_3(\Q(\zeta)^+)\to K_3(\Q(\zeta)^+)/N
}
\]
 Thus, in view of Corollary \ref{cor:suslin},  we can equally regard $\eta_\zeta$ as an element of $B(\Q(\zeta)^+)/N$ and, as such, it can be represented explicitly in terms of the generators of $B(\Q(\zeta)^+)$:  
By \cite[Section 6.4]{hut:cplx13},  if $F$ is a field and if $G\subset \spl{2}{F}$ is a finite cyclic group of order $n$ with generator $t$, the the map $\ho{3}{G}{\Z}\to \ho{3}{\spl{2}{F}}{\Z}\to B(F)$ sends the generator $[\alpha_3(t)]$ to the element
\begin{eqnarray}
\sum_{k=1}^{n-3}\gpb{\frac{t(\infty)-t^{k+1}(\infty)}{t(\infty)-t^{k+2}(\infty)}}+\gpb{1-X}-\gpb{1/X}+\gpb{Z}+\gpb{1/Z}\label{expr}
\end{eqnarray}
where $X=(t(\infty)-y)/(t^{-1}(\infty)-y)$ and $Z=(t(\infty)-y)/(t(\infty)-t(y))$ for some $y\in \projl{F}\setminus G\cdot \infty$. 
(Note: In fact the target of the map in \cite[Section 6.4]{hut:cplx13} is the \emph{refined scissors congruence group}, $RP(F)$,  and the  expression (\ref{expr}) is obtained by composing with the natural homomorphism $RP(F)\to P(F)$ which sends the term $\left\< a\right\>\gpb{x}$ to $\gpb{x}$.)

 Using the fact that $\{ a\}:=\gpb{a}+\gpb{1/a}$ has order dividing $2$ in $B(F)$, and hence vanishes in $B(F)/N$ if $N$ is odd,  and that $\gpb{b}+\gpb{1-b}:=\gpb{0}=-\gpb{\infty}$ is independent of $b\in F^\times$, this can be written
\[
\sum_{k=1}^{n-3}\gpb{\frac{t(\infty)-t^{k+1}(\infty)}{t(\infty)-t^{k+2}(\infty)}}+\gpb{0}=\sum_{k\mod{n}}\gpb{\frac{t(\infty)-t^{k+1}(\infty)}{t(\infty)-t^{k+2}(\infty)}}.
\]

Taking $t=t_\zeta$ as above, we obtain the element
\[
\eta_\zeta=\sum_{k\mod{N}}\gpb{\frac{t(\infty)-t^{k+1}(\infty)}{t(\infty)-t^{k+2}(\infty)}}=\sum_{k\mod{N}}\gpb{1-\left(\frac{\zeta-\zeta^{-1}}{\zeta^k-\zeta^{-k}}\right)^2}
\]
in  $B(\Q(\zeta)^+)/N$ (see \cite[Lemma 5.1]{cgz:bloch} for more details).

Recall that our goal is to show:

\begin{thm}\label{thm:main} Let $N=\ell^m$, $\ell$ an odd prime,  as above.
Then  $c_\zeta(\eta_\zeta)=\zeta$. Equivalently,
\[
\bar{c}_{2,1}(\eta_\zeta)=\zeta\otimes \zeta \in (F_N^\times)/(F_N^\times)^N\otimes \mu_N.
\]
\end{thm}
We begin with the commutative diagram 
\[
\xymatrix{
K_3(\Q(\zeta)^+)/N\ar[r]\ar[d]&K_3(\Q(\zeta)^+;\Z/N)\ar^-{\bar{c}_{2,1}}[r]\ar[d]&H^1(\Q(\zeta)^+,\mu_N^{\otimes 2})\ar[d]\\
K_3(\Q(\zeta))/N\ar[r]&K_3(\Q(\zeta);\Z/N)\ar^-{\bar{c}_{2,1}}[r]&H^1(\Q(\zeta),\mu_N^{\otimes 2})\\
}
\]
Thus, if we let $\bar{\eta}_\zeta$ denote the image of $\eta_\zeta$ in $K_3(\Q(\zeta))/N$ then we must show that $\bar{c}_{2,1}(\bar{\eta}_\zeta)=\zeta\otimes \zeta$. 

\textbf{Caution:} While the map $K_3(\Q(\zeta)^+)/N\to K_3(F_N)/N$ is injective, the map $B(\Q(\zeta)^+)/N\to B(\Q(\zeta))/N$ is the zero map and hence the image of $\bar{\eta}_\zeta$ in $B(\Q(\zeta))/N$ is $0$. 

Let $D:(\Q(\zeta)^+)^\times \to \spl{}{\Q(\zeta)^+}$ be the map $x\mapsto \mathrm{diag}(x^{-1},x)$.  Then $t=t_\zeta$ is conjugate to $D(\zeta)$ in $\spl{2}{\Q(\zeta)}$ and hence the image of 
$[\alpha_3(t)]$ is equal to the image of $[\alpha_3(D(\zeta)]$ in $\ho{3}{\spl{2}{\Q(\zeta)}}{\Z}$.
Thus if we let   $\tau_N:\ho{3}{(D(\mu_N)}{\Z}\to \bar{K}_3(\Q(\zeta))$ denote 
the homomorphism
\[
\xymatrix{H_3(D(\mu_N),\Z)\ar[r]&H_3(\spl{}{\Q(\zeta)},\Z)\ar[r]^(0.65){\bar{h}^{-1}}&K_3(\Q(\zeta))\\},
\]
then $\bar{\eta}=\bar{\eta}_\zeta\in {K}_3(\Q(\zeta))/N$ is the element  (represented by) $\tau_N([\alpha_3(D(\zeta))])$.

\section{The Bott element in $K$-theory}\label{sec:bott}

Let  $R$ be a ring containing a primitive $N$th root of unity $\zeta=\zeta_N$.

The map $\mu_N\to GL_1(R)\to GL(R)$ induces a map of sets 
\[
H_2(\mu_N,\Z/N)=\pi_2(B\mu_N,\Z/N)\to \pi_2(BGL(R),\Z/N)\to \pi_2(BGL(R)^+,\Z/N)=K_2(R,\Z/N). 
\]
This map is a homomorphism when $N$ is odd (see, for example, \cite[IV, Remark 2.5.3]{weibel:kbook}). The \emph{Bott Element} $\beta=\beta(\zeta)\in K_2(R,\Z/N)$ is the image of 
$\tilde{\beta}=\tilde{\beta}(\zeta)\in H_2(\mu_N,\Z/N)$ under this map. 

We note that $\ch{1}{0}(\beta)=\zeta\in H^0(R,\mu_N)=\mu_N$ (see \cite[V, Lemma 11.10.1]{weibel:kbook})). 

\begin{lem}\label{lem:sou}  Let $E$ be a field containing a primitive $N$th root of unity in which $\ell$ is invertible. Identify $H^1(E,\mu_N^{\otimes i})$ with $E^\times/(E^\times)^N\otimes \mu_N^{\otimes (i-1)}$ via the Kummer isomorphism. 
Then the diagram 
\[
\xymatrix{K_3(E,\Z/N)\ar^{\ch{2}{1}}[r]&E^\times/(E^\times)^N\otimes \mu_N\\
K_1(E,\Z/N)\ar^{*\beta}[u]\ar^{\ch{1}{1}}[r]&E^\times/(E^\times)^N\ar^{\otimes\zeta}[u]}
\]
commutes.
\end{lem}

\begin{proof}
Corollary \ref{cor:sou}  gives 
\[
\ch{2}{1}(x*\beta)=-\ch{1}{1}(x)\cup \ch{1}{0}(\beta)=-x\cup \zeta = x\otimes \zeta
\]
for all $x\in K_1(E,\Z/N)$. 
\end{proof}

\begin{cor}\label{cor:zb} Let $E$ be a field containing a primitive $N$th root of unity in which $\ell$ is invertible. \\
 Then $\ch{2}{1}(\zeta *\beta)= \zeta\otimes \zeta \in E^\times/(E^\times)^N\otimes \mu_N\cong H^1(E,\mu_N^{\otimes 2})$.
\end{cor}

We will now show  that $\bar{\eta}_{\zeta}$ maps to $\zeta*\beta$ under $K_3(\Q(\zeta))\to K_3(\Q(\zeta),\Z/N)$  (see Proposition \ref{prop:gamma} below). 
Theorem \ref{thm:main} then follows from Corollary \ref{cor:zb}.

Recall that the maps $D_n:\gl{n}{R}\to \spl{n+1}{R}$, $A\mapsto (A,\det(A)^{-1})$ induce maps on homology 
$D_n:H_\bullet(\gl{n}{R},k)\to H_\bullet(\spl{n+1}{R},k)$ compatible with stabilization (for any commutative ring $R$ and coefficient ring $k$). Taking the limit as $n\to\infty$, we obtain a splitting $D:H_\bullet(\gl{}{R},k)\to H_\bullet(\spl{}{R},k)$ of the 
inclusion $H_\bullet(\spl{}{R},k)\to H_\bullet(\gl{}{R},k)$. 

\begin{lem}\label{lem:hur} For any field $F$ and any odd $N\geq 1$, the Hurewicz homomorphism $h_3:K_3(F;\Z/N)\to H_3(\gl{}{F},\Z/N)$ induces an isomorphism 
$K_3(F;\Z/N)\cong H_3(\spl{}{F},\Z/N)$.
\end{lem}

\begin{proof} Since $\bar{h}_3:\bar{K}_3(F)\to H_3(\spl{}{F},\Z)$ is an isomorphism, it induces an isomorphism $K_3(F)/N\cong H_3(\spl{}{F},\Z)/N$. Thus we have a commutative diagram 
of Bockstein exact sequences where the maps are Hurewicz homomorphims:
\[
\xymatrix{
0\ar[r]&K_3(F)/N\ar[r]\ar^-{\cong}[d]&K_3(F;\Z/N) \ar[r]\ar[d] &K_2(F)[N]\ar[r]\ar^-{\cong}[d] &0\\
0\ar[r] &H_3(\spl{}{F},\Z)/N\ar[r]& H_3(\spl{}{F},\Z/N)\ar[r]&H_2(\spl{}{F},\Z)[N]\ar[r]& 0\\
}
\]
where the right-hand vertical arrow is an isomorphism since the Hurewicz homomorphism induces an isomorphism $K_2(F)\cong H_2(\spl{}{F},\Z)$ for any field $F$.
\end{proof}

Since the map $D:H_3(\gl{}{F},\Z/N)\to H_3(\spl{}{F},\Z/N)$ is a splitting of the inclusion\\
 $H_3(\spl{}{F},\Z/N)\to H_3(\gl{}{F},\Z/N)$, we get:
\begin{cor}\label{cor:hur}  Let $N$ be odd and let $\theta\in H_3(\gl{}{F};\Z/N)$ lie in the image of the Hurewicz homomorphism from $K_3(F;\Z/N)$. Then $\theta=D(\theta)$.
\end{cor}
For any commutative ring $R$, the product structure on $K_\bullet(R;\Z/N)$ is compatible with the Hurewicz homomorphisms $h_n:K_n(R;\Z/N)\to H_n(\gl{}{R},\Z/N)$. The product on these latter groups 
is induced by the tensor product of matrices homomorphism $\gl{n}{R}\times\gl{m}{R}\to \gl{nm}{R}$  (independent, up to conjugation, on choice of ordered bases, and hence compatible with stabilization of the homology groups as $n,m\to\infty$).     In particular, the product on $H_\bullet(R^\times,\Z/N)=H_\bullet(\gl{1}{R},\Z/N)\subset H_\bullet(\gl{}{R},\Z/N)$ is just the Pontryagin product for the abelian group $R^\times$.

By definition of the elements $\zeta\in K_1(\Q(\zeta),\Z/N)$ and $\beta\in K_2(\Q(\zeta),\Z/N\Z)$ we have $h_1(\zeta)=[\zeta]=[\alpha_1(\zeta)]\in H_1(\mu_N,\Z/N)\subset H_1(\gl{}{\Q(\zeta)},\Z/N)$ and 
$h_2(\beta)=\tilde{\beta}=[\alpha_2(\zeta)]\in H_2(\mu_N,\Z/N)\subset H_2(\gl{}{\Q(\zeta)},\Z/N)$. Since the Hurewicz homomorphism $h_n$ respects products, using Corollary \ref{cor:beta} we have:
\begin{cor}\label{cor:hurewicz}
$h_3(\zeta*\beta)=[\zeta]\pont\beta=[\alpha_3(\zeta)]\in H_3(\mu_N,\Z/N)\subset H_3(\gl{}{\Q(\zeta)},\Z/N)$. 
\end{cor}

As remarked above, the following proposition completes the proof of Theorem \ref{thm:main}:
\begin{prop}\label{prop:gamma}  We have  $\bar{\eta}_\zeta= \zeta*\beta$ in $K_3(\Q(\zeta);\Z/N)$.
\end{prop}

\begin{proof} By Lemma \ref{lem:hur}, it is enough to show that $h_3(\bar{\eta})=h_3(\zeta*\beta)$. 
By Corollary \ref{cor:hurewicz}, we have $h_3(\zeta*\beta)=[\alpha_3(\zeta)]$. By Corollary \ref{cor:hur} it follows that $[\alpha_3(\zeta)]=D[\alpha_3(\zeta)]=[\alpha_3(D(\zeta))]$ in 
$H_3(\gl{}{\Q(\zeta)},\Z/N)$. 
On the other hand, by definition of $\eta_\zeta$, we have $h_3(\bar{\eta}_\zeta)=[\alpha_3(D(\zeta))]\in H_3(\gl{}{\Q(\zeta)},\Z/N)$. 
 \end{proof}


\begin{thebibliography}{10}


\bibitem{bass:tate}
{H. Bass, J. Tate}.
\newblock The Milnor Ring of a Global Field.
\newblock Algebraic K-theory, II: "Classical'' algebraic K-theory and connections with arithmetic (Proc. Conf., Seattle, Wash., Battelle Memorial Inst., 1972), pp. 349–446.
\newblock Lecture Notes in Math., Vol. 342, Springer, Berlin, 1973.


\bibitem{brown:coh}
  {Brown, Kenneth S.}.
     \newblock {Cohomology of Groups},
    \newblock{Graduate Texts in Mathematics},{87}, {Springer-Verlag},{New York}, {1982},
   


\bibitem{cgz:bloch}
    {Frank Calegari, Stavros Garoufalidis, Don Zagier},
    \newblock {Bloch Groups, {A}lgebraic {$K$}-theory and {N}ahm's conjecture}.

    \newblock Ann. Sci. École Norm. Sup. 56  (2023), no. 2, 383-426.


\bibitem{hut:cplx13}
Kevin Hutchinson.
\newblock A {B}loch-{W}igner complex for {$\mathrm{SL}\sb 2$}.
\newblock {\em J. K-Theory}, 12(1):15--68, 2013.

\bibitem{levine:k3}
Marc Levine.
\newblock The indecomposable ${K_3}$ of fields.
\newblock Ann. Sci. École Norm. Sup. (4) 22 (1989), no. 2, 255–344.



\bibitem{sah:dupont}
Johan~L. Dupont and Chih~Han Sah.
\newblock Scissors congruences. {II}.
\newblock {\em J. Pure Appl. Algebra}, 25(2):159--195, 1982.
    
\bibitem{soule:etale}
C. Soul\'e
\newblock {$K$}-th\'eorie des anneaux d'entiers de corps de nombres et cohomologie \'etale.
\newblock Invent. Math. 55 (1979), no. 3, 251-295

\bibitem{sus:homgln}
A.~A. Suslin
\newblock Homology of $\mathrm{GL}_n$, characteristic classes and Milnor $K$-theory
\newblock Algebraic K-theory, number theory, geometry and analysis (Bielefeld, 1982), 357–375.
\newblock  Lecture Notes in Math., 1046, Springer, Berlin, 1984.


\bibitem{sus:bloch}
A.~A. Suslin.
\newblock {$K\sb 3$} of a field, and the {B}loch group.
\newblock {\em Trudy Mat. Inst. Steklov.}, 183:180--199, 229, 1990.
\newblock Translated in Proc.\ Steklov Inst.\ Math.\ {\bf 1991}, no.\ 4,
  217--239, Galois theory, rings, algebraic groups and their applications
  (Russian).

\bibitem{weibel:kbook}
Charles A. Weibel
\newblock The $K$-book. An introduction to algebraic $K$-theory
\newblock Graduate Studies in Mathematics 145.
\newblock American Mathematical Society, Providence, RI (2013)
\end{thebibliography}

\end{document}